\documentclass{article}

\newcommand{\be}{\begin{equation}}
\newcommand{\ee}{\end{equation}}
\newcommand{\bea}{\begin{eqnarray}}
\newcommand{\eea}{\end{eqnarray}}
\newcommand{\barray}{\begin{array}}
\newcommand{\earray}{\end{array}}
\newcommand{\pa}{\partial}
\newcommand{\nn}{\nonumber}
\newcommand{\bitem}{\begin{itemize}}
\newcommand{\eitem}{\end{itemize}}
\newtheorem{teo}{Theorem}[section]
\newcommand{\bt}{\begin{teo}}
\newcommand{\et}{\end{teo}}
\newtheorem{Def}{Definition}[section]
\newcommand{\bd}{\begin{Def}}
\newcommand{\ed}{\end{Def}}
\newtheorem{lem}{Lemma}[section]
\newcommand{\bl}{\begin{lem}}
\newcommand{\el}{\end{lem}}
\newtheorem{prop}{Proposition}[section]
\newcommand{\bp}{\begin{prop}}
\newcommand{\ep}{\end{prop}}
\newtheorem{cor}{Corollary}[section]
\newcommand{\bc}{\begin{cor}}
\newcommand{\ec}{\end{cor}}
\newtheorem{ex}{Example}[section]
\newcommand{\bex}{\begin{ex}}
\newcommand{\eex}{\end{ex}}
\newtheorem{rem}{Remark}[section]
\newcommand{\br}{\begin{rem}}
\newcommand{\er}{\end{rem}}

\catcode `\@=11
\@addtoreset{equation}{section}

\begin{document}

\begin{center}
{\Large \textbf{Lax pairs for the equations
describing compatible \\ nonlocal Poisson brackets
of hydrodynamic
type, \\ and integrable reductions of the Lam\'e
equations\footnote{This work was supported by
the Alexander von Humboldt Foundation (Germany),
the Russian Foundation for Basic Research
(grant No. 99--01--00010) and the INTAS
(grant No. 99--1782).}}}
\end{center}

\bigskip
\bigskip

\centerline{\large {O. I. Mokhov}}
\bigskip
\medskip

\section{Introduction. Basic definitions} \label{vved}

In the present work,
the nonlinear equations for
the general nonsingular pairs of compatible nonlocal
Poisson brackets of hydrodynamic type are derived and
the integrability of these equations
by the method of inverse scattering problem
is proved.
For these equations, the Lax pairs with a spectral parameter
are presented. Moreover, we demonstrate the integrability
of the equations for some especially important partial classes
of compatible nonlocal Poisson brackets of hydrodynamic type,
in particular, for the most important
case when one of the compatible Poisson
brackets is local and also for the case when one of the compatible
Poisson brackets is generated by a metric of constant Riemannian
curvature. The case when one of the compatible
Poisson brackets of hydrodynamic type is local and nondegenerate
was studied in detail in our previous paper
\cite{1}, where the corresponding
equations were derived and the integrability of these equations
was announced. This case is very important, since,
as was shown in \cite{1}, any solution
of these equations generates an integrable
bi-Hamiltonian hierarchy of hydrodynamic type by explicit formulae.
Moreover, these equations describe an important class of
integrable reductions of the classical Lam\'e equations.
Accordingly, in the case when one of the compatible
Poisson brackets is generated by a metric of constant Riemannian
curvature, the corresponding equations describe integrable
reductions of the equations for orthogonal curvilinear coordinate systems
in spaces of constant curvature.

\subsection{Local Poisson brackets of hydrodynamic type}

An arbitrary local homogeneous
first-order Poisson bracket, that is,
a Poisson bracket of the form
\be
\{ I,J \} = \int {\delta I \over \delta u^i(x)}
\biggl ( g^{ij}(u(x)) {d \over dx} + b^{ij}_k (u(x))\, u^k_x \biggr )
{\delta J \over \delta u^j(x)} dx,
\label{lok}
\ee
where $u^1,...,u^N$ are local coordinates
on a certain given smooth $N$-dimensional
manifold $M$, is called {\it a local Poisson bracket of
hydrodynamic type} or {a Dubrovin--Novikov bracket} \cite{2}.
Here $u^i(x),\ 1 \leq i \leq N,$ are functions (fields)
of single independent variable $x$, the coefficients
$g^{ij}(u)$ and $b^{ij}_k (u)$ of bracket (\ref{lok})
are smooth functions of local coordinates,
$I[u]$ and $J[u]$ are arbitrary functionals
on the space of fields $u^i(x), \ 1 \leq i \leq N.$
A local bracket (\ref{lok}) is called
{\it nondegenerate} if
$\det (g^{ij} (u)) \not\equiv 0$.

\bt [Dubrovin, Novikov \cite{2}] \label{dn}
If $\det (g^{ij} (u)) \not\equiv 0$, then bracket (\ref{lok})
is a Poisson bracket, that is, it is skew-symmetric and satisfies
the Jacobi identity, if and only if
\bitem
\item [(1)] $g^{ij} (u)$ is an arbitrary flat
pseudo-Riemannian contravariant metric (a metric of
zero Riemannian curvature),

\item [(2)] $b^{ij}_k (u) = - g^{is} (u) \Gamma ^j_{sk} (u),$ where
$\Gamma^j_{sk} (u)$ is the Riemannian connection generated by
the contravariant metric $g^{ij} (u)$
(the Levi--Civita connection).
\eitem
\et

Consequently, for any local nondegenerate Poisson bracket of
hydrodynamic type, there always exist local coordinates
$v^1,...,v^N$ (flat coordinates of the metric $g^{ij}(u)$) in which
all the coefficients of the bracket are constant:
$$
\widetilde g^{ij} (v)
= \eta^{ij} = {\rm \ const}, \ \
\widetilde \Gamma^i_{jk} (v) = 0, \ \
\widetilde b^{ij}_k (v) =0,
$$
that is, the bracket has the form
\be
\{ I,J \} = \int {\delta I \over \delta v^i(x)}
 \eta^{ij} {d \over dx}
{\delta J \over \delta v^j(x)} dx,
\label{(1.6)}
\ee
where $(\eta^{ij})$ is a nondegenerate symmetric constant matrix:
$$
\eta^{ij} = \eta^{ji}, \ \ \eta^{ij} = {\rm const},
\  \ \det \, (\eta^{ij}) \neq 0.
$$

\subsection{Nonlocal Poisson brackets of
hydrodynamic type}

Nonlocal Poisson brackets of hydrodynamic type
(the Mokhov--Ferapontov brackets) were
introduced and studied in the work of the present author
and Ferapontov \cite{3}.
They have the following form:

\be
\{ I,J \} = \int {\delta I \over \delta u^i(x)}
\left ( g^{ij}(u(x)) {d \over dx} + b^{ij}_k (u(x))\, u^k_x
+ K u^i_x \left ( {d \over dx} \right )^{-1} u^j_x \right )
{\delta J \over \delta u^j(x)} dx,
\label{nonl}
\ee
where $K$ is an arbitrary constant.
A bracket of form (\ref{nonl}) is called {\it nondegenerate}
if
$\det (g^{ij} (u)) \not\equiv 0$.

\bt [\cite{3}] \label{mofer}
If $\det (g^{ij} (u)) \not\equiv 0$, then bracket (\ref{nonl})
is a Poisson bracket, that is, it is skew-symmetric and
satisfies the Jacobi identity, if and only if
\bitem
\item [(1)] $g^{ij} (u)$ is an arbitrary pseudo-Riemannian
contravariant metric of
constant Riemannian curvature $K$,

\item [(2)] $b^{ij}_k (u) = - g^{is} (u) \Gamma ^j_{sk} (u),$ where
$\Gamma^j_{sk} (u)$ is the Riemannian connection
generated by the contravariant metric $g^{ij} (u)$
(the Levi--Civita connection).
\eitem
\et

In \cite{4}
Ferapontov introduced and studied more general
nonlocal Poisson brackets of hydrodynamic type
(the Ferapontov brackets), namely,
the Poisson brackets of the form
\bea
&& \
\{ I,J \} = \int   {\delta I \over \delta u^i(x) }
\left ( g^{ij}(u(x)) {d \over dx} +
b^{ij}_k (u(x))\, u^k_x + \right.
\nn\\
&& \
\left. \sum_{\alpha =1}^L
\varepsilon_{\alpha} (w^{\alpha})^i_k (u (x))
u^k_x \left ( {d \over dx} \right )^{-1}
(w^{\alpha})^j_s (u (x)) u^s_x \right )
{\delta J \over \delta u^j(x)} dx,
\ \ \det (g^{ij} (u)) \not\equiv 0,
\label{nonl2}
\eea
where $\varepsilon_{\alpha} = \pm 1,$ $\alpha = 1,...,L.$

\bt [\cite{4}] \label{ferap}
Bracket (\ref{nonl2}) is a Poisson bracket,
that is, it is skew-symmetric and satisfies the
Jacobi identity, if and only if
\bitem
\item [(1)] $b^{ij}_k (u) = - g^{is} (u) \Gamma ^j_{sk} (u),$
where
$\Gamma^j_{sk} (u)$ is the Riemannian connection
generated by the contravariant metric $g^{ij} (u)$
(the Levi--Civita connection),

\item [(2)] the metric $g^{ij} (u)$ and the set of the affinors
$(w^{\alpha})^i_j (u)$ satisfies relations

\be
g_{ik} (u) (w^{\alpha})^k_j (u) =
g_{jk} (u) (w^{\alpha})^k_i (u),
\ \ \ \alpha = 1,...,L,  \label{peter1}
\ee
\be
\nabla_k (w^{\alpha})^i_j (u) =
\nabla_j (w^{\alpha})^i_k (u),
\ \ \ \alpha = 1,...,L,   \label{peter2}
\ee
\be
R^{ij}_{kl} (u) =  \sum_{\alpha =1}^L
\varepsilon_{\alpha}
\left ( (w^{\alpha})^i_l (u) (w^{\alpha})^j_k (u)
- (w^{\alpha})^j_l (u) (w^{\alpha})^i_k (u) \right ). \label{gauss}
\ee
In addition, the family of the affinors $w^{\alpha} (u)$
is commutative: $[w^{\alpha}, w^{\beta}] =0.$
\eitem
\et

Let us write out all the relations on the coefficients
of the nonlocal Poisson bracket
(\ref{nonl2}) in a convenient form for further
repeated use.

\bl
Bracket (\ref{nonl2}) is a Poisson bracket
if and only if its coefficients satisfy the
relations
\be
g^{ij} = g^{ji}, \label{01}
\ee
\be
{\pa g^{ij} \over \pa u^k} = b^{ij}_k + b^{ji}_k, \label{02}
\ee
\be
g^{is} b^{jk}_s = g^{js} b^{ik}_s, \label{03}
\ee
\be
g^{is} (w^{\alpha})^j_s = g^{js} (w^{\alpha})^i_s, \label{04}
\ee
\be
(w^{\alpha})^i_s (w^{\beta})^s_j =
(w^{\beta})^i_s (w^{\alpha})^s_j, \label{05}
\ee
\be
g^{is} g^{jr} {\pa (w^{\alpha})^k_r \over \pa u^s}
- g^{jr} b^{ik}_s (w^{\alpha})^s_r  =
g^{js} g^{ir} {\pa (w^{\alpha})^k_r \over \pa u^s}
- g^{ir} b^{jk}_s (w^{\alpha})^s_r, \label{06}
\ee
\be
g^{is} \left ( {\pa b^{jk}_r \over \pa u^s}
- {\pa b^{jk}_s \over \pa u^r} \right )
+ b^{ik}_s b^{sj}_r - b^{ij}_s b^{sk}_r =
\sum_{\alpha = 1}^L \varepsilon_{\alpha} g^{is}
\left ( (w^{\alpha})^j_s (w^{\alpha})^k_r -
 (w^{\alpha})^j_r (w^{\alpha})^k_s \right ). \label{07}
\ee
\el

\subsection{Compatible Poisson brackets}

In \cite{5} Magri proposed
a bi-Hamiltonian approach to the integration of nonlinear
systems. This approach demonstrated
that integrability is closely related to the
bi-Hamiltonian property, that is, the property
of a system to have two compatible Hamiltonian
representations.

\bd [Magri \cite{5}] \label{dm}
{\rm Two Poisson brackets $\{ \, \cdot \, , \, \cdot \, \}_1$
and $\{ \, \cdot \, , \, \cdot \, \}_2$ are called {\it compatible} if
an arbitrary linear combination of these Poisson brackets
\be
\{ \, \cdot \, , \, \cdot \, \} =
\lambda_1 \, \{ \, \cdot \, , \, \cdot \, \}_1 +
\lambda_2 \, \{ \, \cdot \, , \, \cdot \, \}_2, \label{magri}
\ee
where $\lambda_1$ and $\lambda_2$ are arbitrary constants,
is also a Poisson bracket.
In this case, we shall also say that the brackets
$\{ \, \cdot \, , \, \cdot \, \}_1$ and
$\{ \, \cdot \, , \, \cdot \, \}_2$
form {\it a pencil of Poisson brackets}.}
\ed

 As was shown by Magri in \cite{5},
compatible Poisson brackets generate integrable hierarchies of
systems of differential equations.
In particular, for a system, the bi-Hamiltonian property
generates recurrent relations for the conservation laws of
this system.

\subsection{Compatible pseudo-Riemannian metrics}

Two pseudo-Riemannian contravariant metrics
$g_1^{ij} (u)$ and $g_2^{ij} (u)$ are called {\it compatible}
if for any linear combination of these metrics
$g^{ij} (u) = \lambda_1 g_1^{ij} (u) + \lambda_2 g_2^{ij} (u)$,
where $\lambda_1$ and $\lambda_2$ are arbitrary
constants for which $\det ( g^{ij} (u) ) \not\equiv 0$,
the coefficients of the corresponding
Levi-Civita connections and the components of
the corresponding tensors of Riemannian curvature are related
by the same linear formula:
$\Gamma^{ij}_k (u) = \lambda_1 \Gamma^{ij}_{1, k} (u) +
\lambda_2 \Gamma^{ij}_{2, k} (u)$ and
$R^{ij}_{kl} (u) = \lambda_1 R^{ij}_{1, kl} (u)
+ \lambda_2 R^{ij}_{2, kl} (u)$
(in this case, we shall say also that {\it the metrics
$g^{ij}_1 (u)$ and $g^{ij}_2 (u)$ form a pencil of metrics}) \cite{6},
\cite{7}.
Flat pencils of metrics, that is nothing but
compatible nondegenerate local Poisson brackets
of hydrodynamic type (compatible Dubrovin--Novikov
brackets \cite{2}), were introduced in \cite{8}.
Two pseudo-Riemannian contravariant metrics
$g_1^{ij} (u)$ and $g_2^{ij} (u)$ of constant Riemannian
curvature
$K_1$ and $K_2$ respectively are called {\it compatible}
if any linear combination of these metrics
$g^{ij} (u) = \lambda_1 g_1^{ij} (u) + \lambda_2 g_2^{ij} (u),$
where $\lambda_1$ and $\lambda_2$ are
arbitrary constants for which
$\det ( g^{ij} (u) ) \not\equiv 0$, is
a metric of constant Riemannian curvature
$\lambda_1 K_1 + \lambda_2 K_2$
and the coefficients of the corresponding
Levi-Civita connections are related by the same
linear formula:
$\Gamma^{ij}_k (u) = \lambda_1 \Gamma^{ij}_{1, k} (u) +
\lambda_2 \Gamma^{ij}_{2, k} (u)$ \cite{6}, \cite{7}.
In this case, we shall also say that
{\it the metrics
$g_1^{ij} (u)$ and $g_2^{ij} (u)$ form a pencil of
metrics of constant Riemannian curvature} \cite{6}, \cite{7}.
It is obvious that all these definitions are mutually
consistent, so that if compatible metrics are
metrics of constant Riemannian curvature, then
they form a pencil of metrics of constant Riemannian curvature,
and if compatible metrics are flat, then they form
a flat pencil of metrics. A pair of pseudo-Riemannian
metrics $g_1^{ij} (u)$ and $g_2^{ij} (u)$ is called
{\it nonsingular} if the eigenvalues of this pair
of metrics, that is, the roots of the equation
$\det ( g_1^{ij} (u) -  \lambda g_2^{ij} (u)) =0,$
are distinct (a pencil of metrics which is
formed by a nonsingular pair of metrics is also
called {\it nonsingular}). In \cite{6}, it is proved
that an arbitrary nonsingular pair of metrics is
compatible if and only if there exist local coordinates
$u = (u^1,...,u^N)$ such that both the metrics
are diagonal in these coordinates and have the
following special form (one can consider that
one of the metrics, here
$g^{ij}_2 (u),$ is an arbitrary
diagonal metric):
$g^{ij}_2 (u) = g^i (u) \delta^{ij}$ and
$g^{ij}_1 (u) = f^i (u^i) g^i (u) \delta^{ij},$
where $f^i (u^i),$ $i=1,...,N,$ are functions of
single variable (generally speaking, complex).
It is obvious that the eigenvalues of the
considered pair of metrics are given by the functions
$f^i (u^i),$ $i=1,...,N.$

\section{Compatible metrics and compatible \\
Poisson brackets of hydrodynamic type}

Consider two arbitrary nonlocal Poisson brackets of
hydrodynamic type
\bea
&& \
\{ I,J \}_1 = \int   {\delta I \over \delta u^i(x) }
\left ( g^{ij}_1 (u(x)) {d \over dx} +
b^{ij}_{1, k} (u(x))\, u^k_x + \right.
\nn\\
&& \
\left. \sum_{\alpha =1}^{L_1}
\varepsilon_{1, \alpha} (w^{\alpha}_1)^i_k (u (x))
u^k_x \left ( {d \over dx} \right )^{-1}
(w^{\alpha}_1)^j_s (u (x)) u^s_x \right )
{\delta J \over \delta u^j(x)} dx
\label{nonl2x}
\eea
and
\bea
&& \
\{ I,J \}_2 = \int   {\delta I \over \delta u^i(x) }
\left ( g^{ij}_2 (u(x)) {d \over dx} +
b^{ij}_{2, k} (u(x))\, u^k_x + \right.
\nn\\
&& \
\left. \sum_{\alpha =1}^{L_2}
\varepsilon_{2, \alpha} (w^{\alpha}_2)^i_k (u (x))
u^k_x \left ( {d \over dx} \right )^{-1}
(w^{\alpha}_2)^j_s (u (x)) u^s_x \right )
{\delta J \over \delta u^j(x)} dx.
\label{nonl2y}
\eea

In \cite{6} it was proved that if the Poisson brackets
(\ref{nonl2x}) and (\ref{nonl2y}) are compatible,
then the metrics $g^{ij}_1 (u)$ and $g^{ij}_2 (u)$
of the brackets are also compatible. Moreover,
it was also proved in \cite{6} that if
1) the pair of
metrics $g^{ij}_1 (u)$ and $g^{ij}_2 (u)$
is nonsingular,
2) both the metrics $g^{ij}_1 (u),$ $g^{ij}_2 (u)$ and
the affinors $(w^{\alpha}_1)^i_j (u),$ $(w^{\alpha}_2)^i_j (u)$
can be samultaneously diagonalized in a domain of local coordinates,
then the Poisson brackets are compatible if and only if
the metrics are compatible.
Here we prove in some sense the converse theorem, which is
very important for our method of integrating the equations for
the general nonsingular pairs of arbitrary
compatible nonlocal Poisson brackets of
hydrodynamic type.

\bt   \label{ttt}
If the pair of
metrics $g^{ij}_1 (u)$ and $g^{ij}_2 (u)$
is nonsingular, then the Poisson brackets
$\{I, J \}_1$ and $\{ I, J \}_2$ are compatible if and only
if the metrics are compatible and
both the metrics $g^{ij}_1 (u),$ $g^{ij}_2 (u)$ and
the affinors $(w^{\alpha}_1)^i_j (u),$ $(w^{\alpha}_2)^i_j (u)$
can be samultaneously diagonalized in a domain of local coordinates.
\et

{\it Proof.} It is sufficient to prove here that if
the pair of metrics is nonsingular and the Poisson brackets are
compatible, then
both the metrics $g^{ij}_1 (u),$ $g^{ij}_2 (u)$ and
the affinors $(w^{\alpha}_1)^i_j (u),$ $(w^{\alpha}_2)^i_j (u)$
can be samultaneously diagonalized in a domain of local coordinates.
All the rest was proved in \cite{6}.
First of all, it was proved that in this case
the metrics
$g^{ij}_1 (u)$ and $g^{ij}_2 (u)$ are compatible.
Recall that if the pair of metrics is nonsingular,
then the metrics are compatible if and only if
there exist local coordinates such that the metrics
are diagonal and have the following special form in these
coordinates: $g^{ij}_2 (u) = g^i (u) \delta^{ij}$ and
$g^{ij}_1 (u) = f^i (u^i) g^i (u) \delta^{ij},$
where $f^i (u^i),$ $1 \leq i \leq N,$ are
functions of single variable \cite{6}. The functions $f^i (u^i)$
are the eigenvalues of the pair of metrics $g^{ij}_1 (u)$ and
$g^{ij}_2 (u)$ so that they are distinct by
assumption of the theorem even in the case if they
are constants.
It follows from the compatibility of the
Poisson brackets $\{ I, J \}_1$ and $\{ I, J \}_2$
(it is necessary to consider relation (\ref{04}) for
the pencil $\{ I, J \}_1 + \lambda \{ I, J \}_2$) that
\be
g^{is}_1 (w^{\alpha}_2)^j_s =
g^{js}_1 (w^{\alpha}_2)^i_s, \label{k1}
\ee
\be
g^{is}_2 (w^{\alpha}_1)^j_s =
g^{js}_2 (w^{\alpha}_1)^i_s. \label{k2}
\ee
Besides, from relation (\ref{04}) for
the Poisson brackets $\{ I, J \}_1$ and $\{ I, J \}_2$
we have
\be
g^{is}_1 (w^{\alpha}_1)^j_s =
g^{js}_1 (w^{\alpha}_1)^i_s,  \label{k3}
\ee
\be
g^{is}_2 (w^{\alpha}_2)^j_s =
g^{js}_2 (w^{\alpha}_2)^i_s.   \label{k4}
\ee
So from (\ref{k1}) and (\ref{k4})
in our special local coordinates we get
\be
g^i (w^{\alpha}_2)^j_i =
g^j (w^{\alpha}_2)^i_j,
\ee
\be
f^i (u^i) g^i (w^{\alpha}_2)^j_i =
f^j (u^j) g^j (w^{\alpha}_2)^i_j.
\ee
Thus
\be
(w^{\alpha}_2)^i_j =
{g^i \over g^j} (w^{\alpha}_2)^j_i =
{f^i (u^i) g^i \over f^j (u^j) g^j} (w^{\alpha}_2)^j_i,
\ee
that is,
\be
\left ( 1 - {f^i (u^i) \over f^j (u^j)} \right )
(w^{\alpha}_2)^j_i = 0.
\ee
Consequently, since all the functions $f^i (u^i)$ are distinct, we get
\be
(w^{\alpha}_2)^j_i = 0 \ \ {\rm \ for \  } i \neq j.
\ee
Similarly,
from (\ref{k2}) and (\ref{k3}) we have
\be
(w^{\alpha}_1)^j_i = 0 \ \ {\rm \ for \  } i \neq j.
\ee
Thus theorem \ref{ttt} is proved.

\section{Equations for nonsingular pairs of
compatible \\ nonlocal Poisson brackets of hydrodynamic type}

Consider an arbitrary nonsingular pair of
compatible nonlocal Poisson brackets of hydrodynamic type
(\ref{nonl2x}) and (\ref{nonl2y}), that is, we assume
that the Poisson brackets are compatible and the
pair of metrics $g^{ij}_1 (u)$ and $g^{ij}_2 (u)$ is nonsingular.

\bt
General nonsingular pairs of compatible nonlocal Poisson
brackets of hydrodynamic type are described by the following
consistent integrable nonlinear systems:
\be
{\pa H^{\alpha}_{2, j} \over  \pa u^i} =
\beta_{ij} H^{\alpha}_{2, i}, \ \ \ \ i \neq j, \label{lamx0}
\ee
\be
{\pa \beta_{ij} \over \pa u^k}
=\beta_{ik} \beta_{kj},\ \ \ i\neq j,\ \ i\neq k,\ \ j\neq k, \label{lamx1}
\ee
\be
\epsilon^i_2 {\pa \beta_{ij} \over \pa u^i}+
\epsilon^j_2 {\pa \beta_{ji} \over \pa u^j}+\sum_{s\neq i,\
s\neq j} \epsilon^s_2 \beta_{si} \beta_{sj}
+ \sum_{\alpha = 1}^{L_2} \varepsilon_{2, \alpha}
H^{\alpha}_{2, i} H^{\alpha}_{2, j} = 0,\ \ \ i\neq j. \label{lamx2}
\ee
\bea
&&
\epsilon^i_2 f^i (u^i) {\pa  \beta_{ij}
 \over \pa u^i}+  {1 \over 2} \epsilon^i_2 (f^i)' \beta_{ij} +
\epsilon^j_2 f^j (u^j)
{\pa
\beta_{ji} \over \pa u^j}+
{1 \over 2} \epsilon^j_2 (f^j)' \beta_{ji} + \nn\\
&&
\sum_{s\neq i,\
s\neq j} \epsilon^s_2 f^s (u^s) \beta_{si} \beta_{sj}
+ \sum_{\alpha = 1}^{L_1} \varepsilon_{1, \alpha}
H^{\alpha}_{1, i}
H^{\alpha}_{1, j} = 0 ,\ \ \ i\neq j, \label{lam3}
\eea
\be
{\pa H^{\alpha}_{1, j} \over  \pa u^i} =
\beta_{ij} H^{\alpha}_{1, i}, \ \ \ \ i \neq j, \label{lam00}
\ee
where $f^i (u^i),$ $i=1,...,N,$ are arbitrary given
functions of one variable.
\et

According to theorem \ref{ttt} there exist local coordinates
such that in these coordinates we have
\be
g^{ij}_2 (u) = g^i (u) \delta^{ij}, \ \ \
g^{ij}_1 (u) = f^i (u^i) g^i (u) \delta^{ij},
\ee
\be
(w^{\alpha}_2)^i_j (u) = (w^{\alpha}_2)^i (u) \delta^i_j, \ \ \
(w^{\alpha}_1)^i_j (u) = (w^{\alpha}_1)^i (u) \delta^i_j.
\ee
Moreover, according to the same theorem \ref{ttt}
any pair of Poisson brackets of this form is compatible.
Thus it is sufficient to consider the conditions that
the brackets $\{ I, J \}_1$ and $\{ I, J \}_2$ of
this special form are Poisson brackets.
If the metric $g^{ij} (u)$ and the affinors
$(w^{\alpha})^i_j (u)$ of bracket (\ref{nonl2})
are diagonal, then the conditions (\ref{01})--(\ref{07})
take an especially simple form.
Consider the conditions for the bracket $\{ I, J \}_2$ in
the special local coordinates.

Introduce the standard classical notation
\bea
&&
g^i (u) = {\epsilon^i_2 \over (H_i (u))^2 },\ \ \ d\, s^2 =
\sum_{i=1}^N  \epsilon^i_2
(H_i (u))^2 (d u^i)^2, \\
&&
\beta_{ik} (u) = {1 \over H_i (u)} {\pa H_k \over \pa u^i},\ \ \
i \neq k, \label{vra}
\eea
where $H_i (u)$ are {\it the Lam\'e coefficients}
and
$\beta_{ik} (u)$ are {\it the rotation coefficients},
$\epsilon^i_2 = \pm 1,$ $i = 1,...,N$.
Although, in our case, all the functions
are, generally speaking, complex, we shall use
formulae, which are convenient for using also in
the purely real case.

In the considered ``diagonal'' case,
the conditions that $\{ I, J \}_2$ is a Poisson bracket are
equivalent to the Gauss-Codazzi equations for
submanifolds with flat normal bundle and
holonomic net of curvature lines (see \cite{9}, \cite{10}).
Following \cite{9}, \cite{10} introduce the
functions $H^{\alpha}_{2, i} (u),$
$1 \leq i \leq N,$ $1 \leq \alpha \leq L_2,$ such that
\be
(w^{\alpha}_2)^i (u) = { H^{\alpha}_{2, i} (u) \over
H_i (u) }.
\ee

Then the conditions (\ref{01})--(\ref{07})
for the bracket $\{ I, J \}_2$ take the form (see \cite{9})
\be
{\pa H^{\alpha}_{2, j} \over  \pa u^i} =
\beta_{ij} H^{\alpha}_{2, i}, \ \ \ \ i \neq j, \label{lam0}
\ee
\be
{\pa \beta_{ij} \over \pa u^k}
=\beta_{ik} \beta_{kj},\ \ \ i\neq j,\ \ i\neq k,\ \ j\neq k, \label{lam1}
\ee
\be
\epsilon^i_2 {\pa \beta_{ij} \over \pa u^i}+
\epsilon^j_2 {\pa \beta_{ji} \over \pa u^j}+\sum_{s\neq i,\
s\neq j} \epsilon^s_2 \beta_{si} \beta_{sj}
+ \sum_{\alpha = 1}^{L_2} \varepsilon_{2, \alpha}
H^{\alpha}_{2, i} H^{\alpha}_{2, j} = 0,\ \ \ i\neq j. \label{lam2}
\ee

For the metric
$g^{ij}_1 (u) = f^i (u^i) g^i (u) \delta^{ij}$,
the Lam\'e coefficients and the rotation coefficients
have the form
\bea
&&
\widetilde H_i (u) = {H_i (u)  \over \sqrt {\epsilon^i_1 f^i (u^i)}},
\ \ \ f^i (u^i) g^i (u) = {\epsilon^i_1 \epsilon^i_2 \over
(\widetilde H_i (u))^2},\ \ \ \epsilon^i_1 = \pm 1, \label{redu1}\\
&&
\widetilde \beta_{ik} (u) =
{1 \over \widetilde H_i (u)} {\pa \widetilde H_k \over
\pa u^i} = \nn\\
&&
{\sqrt {\epsilon^i_1 f^i (u^i)} \over
\sqrt {\epsilon^k_1 f^k (u^k)}} \left ( {1 \over H_i (u)} {\pa H_k \over
\pa u^i}  \right ) =  {\sqrt {\epsilon^i_1 f^i (u^i)} \over
\sqrt {\epsilon^k_1 f^k (u^k)}}
\beta_{ik} (u), \  i \neq k. \label{reduction}
\eea
Besides,
introduce the
functions $H^{\alpha}_{1, i} (u),$
$1 \leq i \leq N,$ $1 \leq \alpha \leq L_1,$ such that
\be
(w^{\alpha}_1)^i (u) = {\widetilde H^{\alpha}_{1, i} (u) \over
\widetilde H_i (u) } =
\sqrt {\epsilon^i_1 f^i (u^i)} {\widetilde H^{\alpha}_{1, i} (u) \over
 H_i (u) } = { H^{\alpha}_{1, i} (u) \over
 H_i (u) }.
\ee

Accordingly, equations (\ref{lam1}) are
automatically satisfied also for the rotation coefficients
$\widetilde \beta_{ik} (u)$, equations (\ref{lam0})
become equations (\ref{lam00}), and equations
(\ref{lam2}) for $\widetilde \beta_{ik} (u)$ give equations
(\ref{lam3}).

Note that it is easy to show that equations
(\ref{lamx2}), (\ref{lam3}) for nonsingular pairs of metrics
(that is, all the functions $f^i (u^i)$ must be distinct
also in the case if they are constants) are equivalent
to the following equations (in particular,
it is more convenient to use
these equations for checking the consistency of system
(\ref{lamx0})--(\ref{lam00})):
\bea
&&
{\pa \beta_{ij} \over \pa u^i}
= {1 \over 2} {(f^i (u^i))' \over (f^j (u^j) - f^i (u^i))}
\beta_{ij} + {\epsilon^i_2 \epsilon^j_2 \over 2}
{(f^j (u^j))' \over (f^j (u^j) - f^i (u^i))}
\beta_{ji} - \nn\\
&&
- \sum_{s \neq i, s \neq j} \epsilon^i_2 \epsilon^s_2
{(f^j (u^j) - f^s (u^s)) \over (f^j (u^j) - f^i (u^i))}
\beta_{si} \beta_{sj} + \sum_{\beta = 1}^{L_1} {\epsilon^i_2
\varepsilon_{1, \beta} \over (f^j (u^j) - f^i (u^i))}
 H^{\beta}_{1, i}  H^{\beta}_{1, j} -\nn\\
&&
\sum_{\alpha = 1}^{L_2} {\epsilon^i_2 \varepsilon_{2, \alpha}
f^j (u^j) \over (f^j (u^j) - f^i (u^i))}
H^{\alpha}_{1, i}  H^{\alpha}_{1, j},
 \ \ i \neq j.
\eea

\section{Lax pair for the general nonsingular pair \\
of compatible nonlocal Poisson brackets \\
of hydrodynamic type}

The Lax pair with a spectral parameter for the
system (\ref{lamx0})--(\ref{lam00}) can be derived from
the linear problem for the system (\ref{lamx0})--(\ref{lamx2})
describing all submanifolds with flat normal bundle and
holonomic net of curvature lines. The equations
(\ref{lamx0})--(\ref{lamx2}) are the conditions of consistency
for the following linear system:
\bea
&&
{\pa \varphi_i \over \pa u^k} =
{\sqrt {\epsilon^i_2} \over \sqrt{\epsilon^k_2}} \beta_{ik}
\varphi_k, \ \ \ i \neq k, \label{aaa1}   \\
&&
{\pa \varphi_i \over \pa u^i} = - \sum_{k \neq i}
{\sqrt {\epsilon^k_2} \over \sqrt{\epsilon^i_2}} \beta_{ki}
\varphi_k + \sum_{\alpha = 1}^{L_2}
{\sqrt {\varepsilon_{2, \alpha}} \over \sqrt{\epsilon^i_2}}
H^{\alpha}_{2, i} \psi^{\alpha},   \label{aaa2}  \\
&&
{\pa \psi^{\alpha} \over \pa u^i} = -
{\sqrt {\varepsilon_{2, \alpha}} \over \sqrt{\epsilon^i_2}}
H^{\alpha}_{2, i} \varphi_i.   \label{aaa3}
\eea

The condition that the bracket
$\{ I, J \}_1 + \lambda \{ I, J \}_2$
is a Poisson bracket for any $\lambda$
is equivalent to the system (\ref{lamx0})--(\ref{lamx2})
corresponding to
the metric $(\lambda + f^i (u^i)) g^i (u) \delta^{ij}$ and
the affinors $(w^{\beta}_1)^i_j (u),$
$1 \leq \beta \leq L_1,$ and $\sqrt {\lambda} (w^{\alpha}_2)^i_j (u),$
$1 \leq \alpha \leq L_2$. In this case,
the linear problem (\ref{aaa1})--(\ref{aaa3})
becomes the Lax pair with the spectral parameter $\lambda$
for the general nonsingular pair of arbitrary compatible nonlocal
Poisson brackets of hydrodynamic type:
\bea
&&
{\pa \varphi_i \over \pa u^k} =
{\sqrt {\epsilon^i_2 ( \lambda + f^i)}
\over \sqrt{\epsilon^k_2 (\lambda + f^k)}} \beta_{ik}
\varphi_k, \ \ \ i \neq k,   \label{bbb1} \\
&&
{\pa \varphi_i \over \pa u^i} = - \sum_{k \neq i}
{\sqrt {\epsilon^k_2 (\lambda + f^k)}
\over \sqrt{\epsilon^i_2 (\lambda + f^i)}} \beta_{ki}
\varphi_k + \nn\\
&&
\sum_{\alpha = 1}^{L_2}
{\sqrt {\varepsilon_{2, \alpha} \lambda}
\over \sqrt{\epsilon^i_2 (\lambda + f^i)}}
H^{\alpha}_{2, i} \psi^{\alpha} +
\sum_{\beta = 1}^{L_1}
{\sqrt {\varepsilon_{1, \beta} }
\over \sqrt{\epsilon^i_2 (\lambda + f^i)}}
H^{\beta}_{1, i} \chi^{\beta},   \label{bbb2}   \\
&&
{\pa \psi^{\alpha} \over \pa u^i} = -
{\sqrt {\varepsilon_{2, \alpha} \lambda}
\over \sqrt{\epsilon^i_2 (\lambda + f^i) }}
H^{\alpha}_{2, i} \varphi_i,   \label{bbb3}    \\
&&
{\pa \chi^{\beta} \over \pa u^i} = -
{\sqrt {\varepsilon_{1, \beta} }
\over \sqrt{\epsilon^i_2 (\lambda + f^i) }}
H^{\beta}_{1, i} \varphi_i.  \label{bbb4}
\eea

The Lax pair with the spectral parameter proves
the integrability the system (\ref{lamx0})--(\ref{lam00})
for the general nonsingular pair of compatible nonlocal
Poisson brackets of hydrodynamic type.

If $H^{\alpha}_{1, i} (u) = 0,$
$1 \leq i \leq N,$ $1 \leq \alpha \leq L_1,$
and $H^{\alpha}_{2, i} (u) =0,$ $1 \leq i \leq N,$
$1 \leq \alpha \leq L_2,$
then the system (\ref{lamx0})--(\ref{lam00}) describe
all compatible local Poisson brackets of hydrodynamic type
(compatible Dubrovin--Novikov brackets or flat pencils of metrics).
This system was derived and integrated by
the method of inverse scattering problem in \cite{11}, \cite{12}
with using the Zakharov method of differential reductions \cite{13}.
The Lax pair for this system was demonstrated by Ferapontov
in \cite{14} (here $\varepsilon^i_2 = 1,$ $1 \leq i \leq N$):
\bea
&&
{\pa \varphi_i \over \pa u^j} = \sqrt {(\lambda + f^i)
\over (\lambda + f^j)} \beta_{ij} \varphi_j, \ \ \ i \neq j,
\label{lax1} \\
&&
{\pa \varphi_i \over \pa u^i} = -
\sum_{k \neq i} \sqrt { (\lambda + f^k)
\over (\lambda + f^i)} \beta_{ki} \varphi_k,
\label{lax2}
\eea
where $\lambda$ is a spectral parameter.
The corresponding linear problem for the Lam\'e
equations, that is, equations
(\ref{lam1}), (\ref{lam2}) with $H^{\alpha}_{2, i} (u) =0,$
was well known Darboux yet \cite{15},
see also, for example,
\cite{13}, \cite{16}:
\bea
&&
{\pa \varphi_i \over \pa u^j} =  \beta_{ij} \varphi_j, \ \ \ i \neq j,
\label{lax1a} \\
&&
{\pa \varphi_i \over \pa u^i} = -
\sum_{k \neq i}  \beta_{ki} \varphi_k.
\label{lax2a}
\eea
The condition of consistency for the linear system
(\ref{lax1a}), (\ref{lax2a}) defines the Lam\'e equations.
The Lax pair (\ref{lax1}), (\ref{lax2}) can be
easily derived from the classical linear problem
(\ref{lax1a}), (\ref{lax2a}) for the Lam\'e equations
(see \cite{17}, \cite{18}).

The Lax pair (\ref{lax1}), (\ref{lax2})
is generalized also to the case of arbitrary
nonsingular pencils of metrics of constant Riemannian curvature (see
examples in \cite{14}).
This Lax pair
can be also easily derived from the corresponding linear
problem for the system
describing all the orthogonal curvilinear coordinate
systems in $N$-dimensional spaces of constant curvature $K_2$:
\bea
&&
{\pa \varphi_i \over \pa u^j} = {\sqrt {\varepsilon^i} \over
\sqrt {\varepsilon^j}}
\beta_{ij} \varphi_j, \ \ i \neq j, \label{l1}\\
&&
{\pa \varphi_i \over \pa u^i} = -
\sum_{k \neq i} {\sqrt {\varepsilon^k} \over
\sqrt { \varepsilon^i}}
\beta_{ki} \varphi_k
+ {\sqrt{K_2} \over \sqrt{\varepsilon^i }}
H_i \psi, \label{l2}\\
&&
{\pa \psi \over \pa u^i} = -
 {\sqrt{K_2} \over \sqrt{\varepsilon^i }}
H_i \varphi_i.
\label{l3}
\eea
The condition of consistensy for the linear system
(\ref{l1})--(\ref{l3}) gives the equations for
all orthogonal curvilinear coordinate systems in
$N$-dimensional spaces of constant Riemannian curvature $K_2$.

The corresponding Lax pair
with a spectral parameter for nonsingular pencils
of metrics of constant Riemannian curvature has the form (see \cite{17},
\cite{18}) :
\bea
&&
{\pa \varphi_i \over \pa u^j} = \sqrt {\varepsilon^i (\lambda + f^i)
\over \varepsilon^j (\lambda + f^j)} \beta_{ij} \varphi_j, \ \ \ i \neq j,
\label{lax1b} \\
&&
{\pa \varphi_i \over \pa u^i} = -
\sum_{k \neq i} \sqrt {\varepsilon^k (\lambda + f^k)
\over \varepsilon^i (\lambda + f^i)} \beta_{ki} \varphi_k
+ \sqrt{ \lambda K_2 + K_1 \over \varepsilon^i (\lambda + f^i)}
H_i \psi,
\label{lax2b} \\
&&
{\pa \psi \over \pa u^i} = -
 \sqrt{  \lambda K_2 + K_1 \over \varepsilon^i (\lambda + f^i)}
H_i \varphi_i,
\label{lax3b}
\eea
where $\lambda$ is a spectral parameter.
The condition of consistency for the linear system
(\ref{lax1b})--(\ref{lax3b})
is equivalent to the equations for nonsingular
pencils of metrics of constant Riemannian curvature.

If $H^{\alpha}_{2, i} (u) =0,$ $1 \leq i \leq N,$
$1 \leq \alpha \leq L_1,$
then the corresponding integrable systems
(\ref{lamx0})--(\ref{lam00}) describe
compatible pairs of Poisson brackets of hydrodynamic type
one of which is local. These systems always give integrable reductions
of the classical Lam\'e equations. The corresponding
Lax pairs with a spectral parameter have the form:

\bea
&&
{\pa \varphi_i \over \pa u^k} =
{\sqrt {\epsilon^i_2 ( \lambda + f^i)}
\over \sqrt{\epsilon^k_2 (\lambda + f^k)}} \beta_{ik}
\varphi_k, \ \ \ i \neq k,   \label{ccc1} \\
&&
{\pa \varphi_i \over \pa u^i} = - \sum_{k \neq i}
{\sqrt {\epsilon^k_2 (\lambda + f^k)}
\over \sqrt{\epsilon^i_2 (\lambda + f^i)}} \beta_{ki}
\varphi_k +
\sum_{\beta = 1}^{L_1}
{\sqrt {\varepsilon_{1, \beta} }
\over \sqrt{\epsilon^i_2 (\lambda + f^i)}}
H^{\beta}_{1, i} \chi^{\beta},   \label{ccc2}   \\
&&
{\pa \chi^{\beta} \over \pa u^i} = -
{\sqrt {\varepsilon_{1, \beta} }
\over \sqrt{\epsilon^i_2 (\lambda + f^i) }}
H^{\beta}_{1, i} \varphi_i.  \label{ccc4}
\eea

Similarly, we can construct the Lax pair with
a spectral parameter for any special partial case of
the equations for compatible nonlocal Poisson brackets of hydrodynamic type,
for example, we get integrable reductions of
the equations for orthogonal curvilinear coordinate systems in
spaces of constant Riemannian curvature $K$ (in this case
$H^{\alpha}_{2, i} (u) = \sqrt {K} H_{2, i} (u),$ $L_2 = 1,$ $\alpha = 1$).

\medskip

\begin{flushleft}
Centre for Nonlinear Studies,\\
L.D.Landau Institute for Theoretical Physics, \\
Russian Academy of Sciences\\
e-mail: mokhov@mi.ras.ru; mokhov@landau.ac.ru\\
\end{flushleft}

\end{document}